\documentclass[12pt,a4paper]{article}
\usepackage{amsfonts,amssymb,amsmath,txfonts}
\usepackage[english]{babel}

\DeclareMathOperator{\End}{End}

 \newcommand{\bP}{{\mathbb P}}

\DeclareMathOperator{\E}{E}
\DeclareMathOperator{\GL}{GL}

\DeclareMathOperator{\dis}{\vartriangle}

\begin{document}

\title{Projective Ring Lines and Their Generalisations}

\author{Hans Havlicek\thanks{Email: \texttt{havlicek@geometrie.tuwien.ac.at}}
\\Institute of Discrete Mathematics and Geometry\\ Vienna University of
Technology\\ Wien, Austria}

\maketitle

\begin{abstract}
We give a survey on projective ring lines and some of their substructures which
in turn are more general than a projective line over a ring.

\end{abstract}

\section{Distant graph and connected components}

The \emph{projective line} $\bP(R)$ over any ring $R$ (associative with $1\neq
0$) can be defined in terms of the free left $R$-module $R^2$ as follows
\cite{blunck+he-05}, \cite{herz-95a}: It is the orbit of a starter point
$R(1,0)$ under the action of the general linear group $\GL_2(R)$ on $R^2$. A
basic notion on $\bP(R)$ is its \emph{distant relation}: Two points are called
distant (in symbols: $\dis$) if they can be represented by the elements of a
two-element basis of $R^2$. The \emph{distant graph} $(\bP(R),\dis)$ has as
vertices the points of $\bP(R)$ and as edges the pairs of distant points. The
distant graph is connected precisely when $\GL_2(R)$ is generated by the
\emph{elementary linear group} $\E_2(R)$, i.e., the subgroup of $\GL_2(R)$
which is generated by elementary transvections, together with the set of all
invertible diagonal matrices \cite{blunck+h-01a}. The orbit of $R(1,0)$ under
$\E_2(R)$ is a connected component of the distant graph. It admits a
parametrisation in terms of infinitely many formulas \cite{blunck+h-01a},
\cite{blunck+h-03b}. The situation is less intricate for a ring $R$ of
\emph{stable rank} $2$ (see \cite{chen-11a}, \cite{veld-85}, or
\cite{veld-95}), as it gives rise to a connected distant graph with diameter
$\leq 2$. The above-mentioned parametrisation turns into \emph{Bartolone's
parametrisation} \cite{bart-89} of $\bP(R)$, namely
\begin{displaymath}
    \bP(R) = \{ R(t_2t_1-1,t_2)\mid t_1,t_2\in R \}
    \mbox{~~~~~($R$ of stable rank $2$)} .
\end{displaymath}
Refer to the seminal paper of P.~M.~Cohn \cite{cohn-66} for the algebraic
background, and to the work of A.~Blunck \cite{blunck-97a}, \cite{blunck-02a}
for orbits of the point $R(1,0)$ under other subgroups of $\GL_2(R)$. 

\section{Chain Geometries, subspaces and Jordan Systems}

Let $R$ be an algebra over a commutative field $K$; by identifying $K$ with
$K\cdot 1_R$ the projective line $\bP(K)$ is embedded in $\bP(R)$. For $R\neq
K$ the projective line $\bP(R)$ can be considered as the point set of the
\emph{chain geometry} $\Sigma(K,R)$; the $\GL_2(R)$ orbit of $\bP(K)$ is the
set of \emph{chains} \cite{blunck+he-05}, \cite{herz-95a}. The geometries of
M\"{o}bius, Minkowski and Laguerre are well known examples of chain geometries
\cite{benz-73}. A crucial property is that any three mutually distinct points
are on a unique chain. The chain geometry $\Sigma(K,R)$ may be viewed as a
refinement of the distant graph, since two points of $\bP(R)$ are distant if,
and only if, they are on a common chain. There are cases though, when the word
``refinement'' is inappropriate in its strict sense: Let $R=\End_F(V)$ be the
endomorphism ring of a vector space $V$ over a (not necessarily commutative)
field $F$ and let $K$ denote the \emph{centre} of $F$. Then the $K$-chains of
$\bP(R)$ can be defined solely in terms of the distant graph $(\bP(R),\dis)$
\cite{blunck+h-12z}.

Each chain geometry $\Sigma(K,R)$ is a \emph{chain space}; see
\cite{blunck+he-05}, where also the precise definition of \emph{subspaces} of a
chain space is given. The algebraic description of subspaces of $\Sigma(K,R)$
is due to A.~Herzer \cite{herz-92b} and H.-J.~Kroll \cite{kroll-91a},
\cite{kroll-92b}, \cite{kroll-92a}. It is based on the following notions: A
\emph{Jordan system} is a $K$-subspace of $R$ satisfying two extra conditions:
(i) $1\in J$; (ii) If $b\in J$ has an inverse in $R$ then $b^{-1}\in J$. (See
\cite{loos-75} for relations with \emph{Jordan algebras} and \emph{Jordan
pairs} and compare with \cite{gold-06a}, \cite{matt-07a}.) A Jordan system $J$
is called \emph{strong} if it satisfies a (somewhat technical) condition which
guarantees the existence of ``many'' invertible elements in $J$. Strong Jordan
systems are closed under \emph{triple multiplication}, i.~e., $ xyx\in J$ for
all $x,y\in J$. The \emph{projective line $\bP(J)$ over a strong Jordan system}
$J\subset R$ is defined by restricting the \emph{parameters} $t_1,t_2$ to $J$
in Bartolone's parametrisation. We wish to emphasise that in general a point of
$\bP(J)$ cannot be written as $R(a,b)$ with $a,b\in J$, unless $J$ is even a
subalgebra of $R$. The main theorem about subspaces is as follows: If $R$ is a
strong algebra then any connected subspace of $\Sigma(K,R)$ is projectively
equivalent to a projective line over a strong Jordan system of $R$.

Projective lines over strong Jordan systems admit many applications: For
example, one may use them to describe subsets of Grassmannians which are closed
under reguli \cite{herz-92b} or chain spaces on quadrics \cite{blunck-97}. See
also \cite{blunck-94}, \cite{herz-08a}, \cite{herz-10a}, \cite{herz-11a}, and
the numerous examples given in \cite{blunck+he-05}.
\par
Finally, let us mention one of the many questions that remain: \emph{Is it
possible to replace the strongness condition for Jordan systems by closedness
under triple multiplication without affecting the known results?} A partial
affirmative answer was given in \cite{blunck+h-10a} for case when $R$ is the
ring of $n\times n$ matrices over a field $F$ with an involution $\sigma$ and
$J$ is the (not necessarily strong) Jordan system of $\sigma$-Hermitian
matrices. The proof is based upon the verification that the projective line
over this $J$ is, up to some notational differences, nothing but the point set
of a \emph{dual polar space} \cite{cameron-82a} or, in the terminology of
\cite{wan-96}, the point set of a \emph{projective space of $\sigma$-Hermitian
matrices}.
\par
A wealth of further references can be found in \cite{benz-73},
\cite{blunck+he-05}, \cite{havl-07a}, \cite{herz-95a}, \cite{huanglp-06a},
\cite{pank-10a}, \cite{veld-95}, and \cite{wan-96}.
Refer to \cite{brehm-08}, \cite{brehm+g+s-95}, \cite{faure-04a},
\cite{havl+m+p-11a}, \cite{havl+k+o-12z}, \cite{havlicek+saniga-09a}, and
\cite{lash-97} for deviating definitions of projective lines which we cannot
present here.


\begin{thebibliography}{10}

\bibitem{bart-89} C.~Bartolone.
\newblock Jordan homomorphisms, chain geometries and the fundamental theorem.
\newblock {\em Abh.\ Math.\ Sem.\ Univ.\ Hamburg}, 59:93--99, 1989.

\bibitem{benz-73} W.~Benz.
\newblock {\em {V}orlesungen \"{u}ber {G}eometrie der {A}lgebren}.
\newblock Springer, Berlin, 1973.

\bibitem{blunck-94} A.~Blunck.
\newblock Chain spaces over {J}ordan systems.
\newblock {\em Abh.\ Math.\ Sem.\ Univ.\ Hamburg}, 64:33--49, 1994.

\bibitem{blunck-97} A.~Blunck.
\newblock Chain spaces via {C}lifford algebras.
\newblock {\em Monatsh.\ Math.}, 123:98--107, 1997.

\bibitem{blunck-97a} A.~Blunck.
\newblock {\em Geometries for Certain Linear Groups over Rings --- Construction
  and Coordinatization}.
\newblock Habilitationsschrift, Technische Universit\"{a}t Darmstadt, 1997.

\bibitem{blunck-02a} A.~Blunck.
\newblock Projective groups over rings.
\newblock {\em J.\ Algebra}, 249:266--290, 2002.

\bibitem{blunck+h-01a} A.~Blunck and H.~Havlicek.
\newblock The connected components of the projective line over a ring.
\newblock {\em Adv.\ Geom.}, 1:107--117, 2001.

\bibitem{blunck+h-03b} A.~Blunck and H.~Havlicek.
\newblock Jordan homomorphisms and harmonic mappings.
\newblock {\em Monatsh.\ Math.}, 139:111--127, 2003.

\bibitem{blunck+h-10a} A.~Blunck and H.~Havlicek.
\newblock Projective lines over {J}ordan systems and geometry of {H}ermitian
  matrices.
\newblock {\em Linear Algebra Appl.}, 433:672--680, 2010.

\bibitem{blunck+h-12z} A.~Blunck and H.~Havlicek.
\newblock Geometric structures on finite- and infinite-dimensional
  {G}rassmannians.
\newblock {\em Beitr. Algebra Geom.}, to appear.

\bibitem{blunck+he-05} A.~Blunck and A.~Herzer.
\newblock {\em Kettengeometrien -- {E}ine {E}inf\"{u}hrung}.
\newblock Shaker Verlag, Aachen, 2005.

\bibitem{brehm-08} U.~Brehm.
\newblock Algebraic representation of mappings between submodule lattices.
\newblock {\em J. Math. Sci. (N. Y.)}, 153(4):454--480, 2008.

\bibitem{brehm+g+s-95} U.\ Brehm, M.\ Greferath, and S.~E. Schmidt.
\newblock Projective geometry on modular lattices.
\newblock In F.\ Buekenhout, editor, {\em Handbook of Incidence Geometry},
  pages 1115--1142. Elsevier, Amsterdam, 1995.

\bibitem{cameron-82a} P.~J. Cameron.
\newblock Dual polar spaces.
\newblock {\em Geom. Dedicata}, 12(1):75--85, 1982.

\bibitem{chen-11a} H.~Chen.
\newblock {\em Rings Related to Stable Range Conditions}, volume~11 of {\em
  Series in Algebra}.
\newblock World Scientific, Singapore, 2011.

\bibitem{cohn-66} P.~M. Cohn.
\newblock On the structure of the {$\text{GL}_2$} of a ring.
\newblock {\em Inst.\ Hautes Etudes Sci.\ Publ.\ Math.}, 30:365--413, 1966.

\bibitem{faure-04a} C.-A. Faure.
\newblock Morphisms of projective spaces over rings.
\newblock {\em Adv. Geom.}, 4(1):19--31, 2004.

\bibitem{gold-06a} D.~Goldstein, R.~M. Guralnick, L.~Small, and E.~Zelmanov.
\newblock Inversion invariant additive subgroups of division rings.
\newblock {\em Pacific J. Math.}, 227(2):287--294, 2006.

\bibitem{havl-07a} H.~Havlicek.
\newblock From pentacyclic coordinates to chain geometries, and back.
\newblock {\em Mitt.\ Math.\ Ges.\ Hamburg}, 26:75--94, 2007.

\bibitem{havl+m+p-11a} H.~Havlicek, A.~Matra{\'s}, and M.~Pankov.
\newblock Geometry of free cyclic submodules over ternions.
\newblock {\em Abh. Math. Semin. Univ. Hambg.}, 81(2):237--249, 2011.

\bibitem{havl+k+o-12z} H.~Havlicek, J.~Kosiorek, and B.~Odehnal.
\newblock A point model for the free cyclic submodules over ternions.
\newblock {\em Results Math.}, to appear.

\bibitem{havlicek+saniga-09a} H.~Havlicek and M.~Saniga.
\newblock Vectors, cyclic submodules, and projective spaces linked with
  ternions.
\newblock {\em J. Geom.}, 92(1-2):79--90, 2009.

\bibitem{herz-92b} A.~Herzer.
\newblock On sets of subspaces closed under reguli.
\newblock {\em Geom.\ Dedicata}, 41:89--99, 1992.

\bibitem{herz-95a} A.~Herzer.
\newblock Chain geometries.
\newblock In F.\ Buekenhout, editor, {\em Handbook of Incidence Geometry},
  pages 781--842. Elsevier, Amsterdam, 1995.

\bibitem{herz-08a} A.~Herzer.
\newblock Konstruktion von {J}ordansystemen.
\newblock {\em Mitt. Math. Ges. Hamburg}, 27:203--210, 2008.

\bibitem{herz-10a} A.~Herzer.
\newblock Die kleine projektive {G}ruppe zu einem {J}ordansystem.
\newblock {\em Mitt. Math. Ges. Hamburg}, 29:157--168, 2010.

\bibitem{herz-11a} A.~Herzer.
\newblock Korrektur und {E}rg\"anzung zum {A}rtikel \emph{{D}ie kleine
  projektive {G}ruppe zu einem {J}ordansystem} in {M}itt. {M}ath. {G}es.
  {H}amburg 29, {A}rmin {H}erzer.
\newblock {\em Mitt. Math. Ges. Hamburg}, 30:15--17, 2011.

\bibitem{huanglp-06a} L.-P. Huang.
\newblock {\em Geometry of Matrices over Ring}.
\newblock Science Press, Beijing, 2006.

\bibitem{kroll-91a} H.-J. Kroll.
\newblock Unterr\"{a}ume von {K}ettengeometrien und {K}ettengeometrien mit
  {Q}uadrikenmodell.
\newblock {\em Results Math.}, 19:327--334, 1991.

\bibitem{kroll-92b} H.-J. Kroll.
\newblock Unterr\"{a}ume von {K}ettengeometrien.
\newblock In N.~K. Stephanidis, editor, {\em Proceedings of the 3rd Congress of
  Geometry (Thessaloniki, 1991)}, pages 245--247, Thessaloniki, 1992. Aristotle
  Univ.

\bibitem{kroll-92a} H.-J. Kroll.
\newblock Zur {D}arstellung der {U}nterr\"{a}ume von {K}ettengeometrien.
\newblock {\em Geom.\ Dedicata}, 43:59--64, 1992.

\bibitem{lash-97} A.~Lashkhi.
\newblock Harmonic maps over rings.
\newblock {\em Georgian Math.\ J.}, 4:41--64, 1997.

\bibitem{loos-75} O.~Loos.
\newblock {\em Jordan Pairs}, volume 460 of {\em Lecture Notes in Mathematics}.
\newblock Springer, Berlin, 1975.

\bibitem{matt-07a} S.~Mattarei.
\newblock Inverse-closed additive subgroups of fields.
\newblock {\em Israel J. Math.}, 159:343--347, 2007.

\bibitem{pank-10a} M.~Pankov.
\newblock {\em {G}rassmannians of Classical Buildings}, volume~2 of {\em
  Algebra and Discrete Mathematics}.
\newblock World Scientific, Singapore, 2010.

\bibitem{veld-85} F.~D. Veldkamp.
\newblock Projective ring planes and their homomorphisms.
\newblock In R.~Kaya, P.~Plaumann, and K.~Strambach, editors, {\em Rings and
  Geometry}, pages 289--350. D.\ Reidel, Dordrecht, 1985.

\bibitem{veld-95} F.~D. Veldkamp.
\newblock Geometry over rings.
\newblock In F.\ Buekenhout, editor, {\em Handbook of Incidence Geometry},
  pages 1033--1084. Elsevier, Amsterdam, 1995.

\bibitem{wan-96} Z.-X. Wan.
\newblock {\em Geometry of Matrices}.
\newblock World Scientific, Singapore, 1996.

\end{thebibliography}

\end{document}